\documentclass[12pt]{article}
\usepackage{amssymb}

\newtheorem{theorem}{Theorem}
\newtheorem{corollary}{Corollary}

\title{On the image of convolutions along an arithmetic progression}

\author{Ernie Croot and Chi-Nuo Lee}

\begin{document}

\maketitle

\begin{abstract}  We consider the question of determining
the structure of the set of all $d$-dimensional vectors of the form $N^{-1}(1_A*1_{-A}(x_1), ..., 1_A*1_{-A}(x_d))$
for $A \subseteq \{1,...,N\}$, and also the set of
all $(2N+1)^{-1}(1_B*1_B(x_1), ..., 1_B*1_B(x_d))$, for 
$B \subseteq \{-N, -N+1, ..., 0, 1, ..., N\}$, where
$x_1,...,x_d$ are fixed positive integers (we let 
$N \to \infty$).  
Using an elementary method related to the Birkhoff-von Neumann theorem on decompositions of doubly-stochastic
matrices we show that both the above two sets of vectors
roughly form polytopes; and of particular interest is 
the question of bounding the number of corner vertices,
as well as understand their structure.
\end{abstract}

\section{Introduction}

Fix an additive group ${\mathbb G}$, and suppose $A, B$ are finite
subsets of ${\mathbb G}$.  Understanding the structure of sumsets $A+B := \{a+b\ :\ a\in A, b\in B\}$ is an important theme in additive combinatorics and additive number theory.  And one way this is often done is through studying convolutions 
$$
1_A*1_B(x)\ =\ \sum_{u+v = x \atop u,v \in {\mathbb G}} 1_A(u)
1_B(v)\ =\ \#\{a \in A,\ b \in B\ :\ a+b =x\}.
$$
(Although there may be infinitely many $u,v$ in the case
where ${\mathbb G}$ is infinite, all but a finite number of the 
terms will be $0$.)

A basic question one could ask about the distribution of these convolutions is the following.  Suppose we take 
${\mathbb G} = {\mathbb Z}$, and suppose we fix some
distinct points $x_1, x_2, ..., x_d \in {\mathbb Z}$.  
What can one say about the image of the mapping?
$$
f\ :\ A\subseteq \{1,2,...,N\}\ \longrightarrow\ 
(1_A*1_A(x_1),\ 1_A*1_A(x_2),\ ...,\ 1_A*1_A(x_d)).
$$

Why might we care about this problem, besides the fact that 
it's a very natural one to pose?  To answer this, suppose
we knew the values of $1_A*1_A(x)$ for all $x$ in a subset
of the $\{x_1,...,x_d\}$.  Could we maybe then conclude 
something about the convolution at the remaining $x_i$'s?  
One could imagine a strong enough solution to this kind of 
problem could have some applications in solving other additive
problems.  For example, suppose we knew that $A$ was not too small, say 
$A \subseteq [-N, N]$, $|A| > N d^{-\varepsilon}$.
For which $x_1,...,x_d$ would knowing that 
$1_A*1_A(x_1), ..., 1_A*1_A(x_{d-1})$ are all near $0$ imply that $1_A*1_A(x_d)$ is near $0$?  A good enough
answer to this might help better understand when $A+A$
contains long arithmetic progressions, a well-studied 
problem \cite{bourgain, croot1, croot2, green}.
\bigskip

To address the original question of understanding
the image of $f$, we define the sets of normalized convolutions given integers $N \geq 1$ as follows.
$$
S_N(x_1,...,x_d)\ :=\ \{N^{-1} (1_A*1_{-A}(x_1),\ ...,\ 1_A*1_{-A}(x_d))\ :\ A \subseteq \{1,2,...,N\}\}
$$
and
\begin{eqnarray*}
T_N(x_1,...,x_d)\ :=\ &&\{(2N+1)^{-1}(1_A*1_A(x_1),\ ...,\ 1_A*1_A(x_dd))\ :\\
&&\hskip0.5in A \subseteq \{-N, ..., -1, 0, 1,...,N\}\}.
\end{eqnarray*}
Note that 
$$
S_N(x_1,...,x_d),\ T_N(x_1,...,x_d)\ \subseteq\ [0,1]^d.
$$
We will first be working with the case $x_i=i$, $i=1,...,d$,
which motivates why we chose to define $T_N(x_1,...,x_d)$
in terms of subsets $A \subset [-N,N]$ instead of $[1,N]$.
If $A$ consisted of positive integers 
then the convolutions $1_A*1_A(i)$ would be very small, making the problem less interesting.
\bigskip

The main focus of this work is to show that
$S$ and $T$ are approximately polytopes, and to develop 
descriptions of the corner vertices in this polytope, as well as bounds on their number.  To accomplish this we
use a variant of the Birkhoff-von-Neumann Theorem 
\cite{birkhoff},
as discussed in section \ref{weightedgraphs}.

We will also be working with unnormalized $\ell^\infty$
norms of vectors, which we define as
$$
\|(u_1,...,u_d)\|_\infty\ =\ \max_{i=1,...,d} |u_i|.
$$
We will the slightly non-standard notation
$$
{\mathbb Z}_M\ =\ {\rm the\ cyclic\ group\ of\ order\ M
\ under\ addition}.
$$
And for a set of vectors $\vec v_1,...,\vec v_d$ in
a real vector space we define
\begin{eqnarray*}
&&{\rm convex-hull}(\{v_1,...,v_d\})\ =\ \\
&&\hskip0.5in \{ 
\lambda_1 \vec v_1 + \cdots + \lambda_d \vec v_d\ :\ 
\lambda_1,...,\lambda_d \geq 0,\ \lambda_1 + \cdots +
\lambda_d\ =\ 1\}.
\end{eqnarray*}
Our main theorem is as follows.

\begin{theorem} \label{main_theorem}
Let $S = S_N(1,2,3,...,d)$ and $T = T_N(1,2,3,..., d)$; 
that is, we are working with the case where $x_i = i$, $i=1,...,d$.
\bigskip

\begin{itemize} 

\item {\bf Shape of $S$ and $T$:}  
There exist sequences of points 
$$
\vec y_1,...,\vec y_m, \vec z_1, ..., \vec z_{m'}\ \in\ [0,1]^d,
$$
such that if 
$$
H\ :=\ {\rm convex}-{\rm hull}(\{\vec y_1, ..., \vec y_m\}),
\ H'\ :=\ {\rm convex -} {\rm hull}(\{\vec z_1, ..., \vec z_{m'}\}),
$$ 
then for any $\vec s \in S$ and $\vec t \in T$
there exist $\vec h \in H$ and $\vec h' \in H'$ 
such that
\begin{equation}\label{AC}
\| \vec s - \vec h \|_\infty\ =\ O(2^dm / N),\ 
\| \vec t - \vec h' \|_\infty\ =\ O(2^dm / N).
\end{equation}
And, conversely, for any $\vec h \in H$ and $\vec h' \in H'$ there exist 
$\vec s\in S$ and $\vec t \in T$ 
so that (\ref{AC}) holds.

\item {\bf Number of Corners:}  We will show that 
the polytopes $H$ and $H'$ have at most
$$
m\ \leq\ 2^{d(d+1)},\ m'\ \leq\ 4^{d(d+1)}
$$
corners, respectively.  

\item {\bf Rational coordinates and subsets of 
cyclic groups:}
Furthermore, we will show that the $\vec y_i$s and $z_j$s
have the following property:  for each $i=1,...,m$ and  
$j=1,...,m'$, there exist positive integers $M = M(i) \leq 2^d$, $M' = M'(j) \leq 4^d$, and subsets 
$B_i \subseteq {\mathbb Z}_M$ and $B'_j \subseteq 
{\mathbb Z}_{M'}$, such that
$$
\vec y_i\ =\ (a_{i,1}/M, ..., a_{i,d} / M),\ 
 \vec z_j\ =\ (b_{j,1}/M', ..., b_{j,d}/M')\ \in\ [0,1]^d,
$$
where for $i,j,h=1,...,d$,
$$
a_{i,h}\ =\ 1_{B_i}*1_{-B_i}(h),\ b_{j,h}\ =\ 1_{B'_j}*1_{B'_j}(h).
$$
\end{itemize}
\end{theorem}

We have the following immediate corollary.

\begin{corollary} \label{corollary1} For arbitrary positive integers 
$x_1 < x_2 < \cdots < x_d$ we have that $S_N(x_1,...,x_d)$ and
$T_N(x_1,...,x_d)$ are approximately polytopes
as in Theorem \ref{main_theorem}; however, the number of
corners $k$ and $k'$ in this case will be bounded as follows:
$$
k\ \leq\ 2^{x_d(x_d+1)},\ k'\ \leq\ 4^{x_d(x_d+1)}.
$$
\end{corollary}

Now we discuss how to prove the corollary.
We first consider the case $S_N(x_1,...,x_d)$ and 
bounds on $k$.  
To prove the corollary in this case, 
first let $\vec v \in S_N(x_1,...,x_d)$.
Note that this $\vec v$ has an associated subset
$A \subseteq \{1,...,N\}$ so that the $i$th coordinate of
$\vec v$ equals $N^{-1} 1_A*1_{-A}(x_i)$.  This $i$th 
coordinate is, in turn, the $x_i$th coordinate of some 
vector $\vec w \in S_N(1,2,...,x_d)$ of dimension
$x_d$.  This motivates defining the projection operator
$$
\rho\ :\ {\mathbb R}^{x_d}\ \to\ {\mathbb R}^d,
$$
where 
$$
\rho((z_1, z_2, ..., z_{x_d}))\ =\ (z_{x_1}, z_{x_2}, ..., 
z_{x_d}).
$$
We note that $\rho(S_N(1,2,...,x_d)) = S_N(x_1,...,x_d)$,
and if $\vec c_1, ..., \vec c_m$ are corners whose
convex hull contains all of $S_N(1,2,...,x_d)$, then the
convex hull of $\rho(\vec c_1), ..., \rho(\vec c_m)$ 
will contain all of $S_N(x_1,...,x_d)$.  Among 
$\rho(\vec c_1), ..., \rho(\vec c_m)$ will be a minimal
set of $k$ vectors whose convex hull contains $S_N(x_1,...,x_d)$.  So we have
$$
k\ \leq\ m\ \leq\ 2^{x_d(x_d+1)}.
$$
by Theorem \ref{main_theorem}.  The proof for the bounds on
$k'$ are similar, giving
$$
k'\ \leq\ m'\ \leq\ 4^{x_d(x_d+1)}.
$$

\subsection{Discussion and the special case $d=2$}

To get a feel for what this theorem is saying, we consider the case $d=2$ and only the convolutions $1_A*1_{-A}$.  
In this case, it turns out that the set $H$ that encloses
$S$ is the convex hull of the following points in ${\mathbb R}^2$:
\begin{equation}\label{thepointshere}
(0,0),\ (1,1),\ (1/4, 0),\ (0,1/2).
\end{equation}

Now we will see why this is the case:  the theorem claims that the corners are determined by sets
$B_i \subseteq \{0,1,...,M_i-1\}$, for some integer 
$M_i \geq 1$.  Based on the proof, it will turn out that 
the possible $M_i$
here are the cycle lengths in a de Bruijn graph for binary
strings of length $d=2$; and each corner (and set $B_i$) will correspond to one of these cycles.  This graph has 
$4$ vertices, $8$ edges, and $6$ cycles.  The possible cycle lengths for this de Bruijn graph are $1, 2, 3,$ and $4$, and nothing else (the cycle of length $1$ is via the loops from the vertex for the strings $00$ and $11$ to themselves.)  The $6$ possible cycles (including loops), along with the associated $M$, sets $B_i$, and $M^{-1}(1_{B_i}*1_{-B_i}(1), 
1_{B_i}*1_{-B_i}(2))$ they correspond to, respectively, 
are:  
\begin{eqnarray*}
&& 00\ \to\ 00,\ M = 1,\ B = \emptyset,\ (0,0)\\
&& 11\ \to\ 11,\ M = 1,\ B = \{0\},\ (1,1) \\
&&\ 01\ \to\ 10\ \to\ 01,\ M = 2,\ B = \{0\},\ (0,1/2) \\
&& 00\ \to\ 01\ \to\ 10\ \to\ 00,\ M = 3,\ 
B = \{0\},\ (0,0)\\
&& 11\ \to\ 10\ \to\ 01\ \to\ 11,\ M = 3,\ B = \{0,1\},\ 
(1/3, 1/3)\\
&& 00\ \to\ 01\ \to\ 11\ \to\ 10\ \to\ 00,\ M = 4,\ 
B = \{0,1\},\ (1/4, 0).
\end{eqnarray*}

It turns out that the point $(1/3, 1/3)$ is contained within the convex hull of the points (\ref{thepointshere});
and this is minimal, as we cannot make the list any smaller
(each of the $4$ remaining points cannot be enclosed in
the convex hull of the remaining $3$).  Thus, 
(\ref{thepointshere}) are exactly the corners of $H$.

We note that $2^{d(d+1)} = 64$ is much larger than $4$, the number of corners we use.  Thus, already with $d=2$ we see that this bound is far from being tight.

\subsection{Some unsolved problems and discussion} \label{unsolved_section}

We finish the introduction by introducing some problems that we were not able to solve.

\begin{enumerate}

\item A theorem similar to Theorem \ref{main_theorem}
might be provable using the circle method.  One could 
imagine perhaps the ``corners" of the polytope might 
fall out of some ``major arc" calculations.  It might
be worth exploring whether the {\it reverse} is true,
however:  are there any interesting problems traditionally
solved using the circle method that can be (possibly more
naturally) solved using cycles in graphs along with 
Birkhoff-von Neumann decompositions as we do here?  

\item Determine the best possible bounds for $m$ and $m'$, the number of corners needed for the polytopes in the above theorem.  Perhaps the true upper bound is something like $m, m' \leq c\cdot 2^d$.

\item What is the best upper bound on $k$ and $k'$ in 
Corollary \ref{corollary1} one can prove?  Can one prove
$k, k'$ are bounded from above by a function of $d$,
perhaps $k, k' \leq c 2^d$ (for some $c>0$)?  
One case where it might be true is when $x_i = i M$, $i=1,...,d$.  For example,
if $N > M$ is a prime number then we will have that for any
set $A \subseteq {\mathbb Z}_N$ (switching temporarily to 
when the ambient group is ${\mathbb Z}_N$) 
and $B \equiv t \cdot A
\pmod{N}$ with $t \equiv M^{-1} \pmod{N}$ 
($B$ is the set of dilates of $A$ by the multiplicative 
inverse of $M$), then
$$
1_A*1_{-A}(x_i)\ =\ 1_{t A}*1_{-tA}(t x_i)\ =\ 
1_B*1_{-B}(i). 
$$
And for $M$ much smaller than $N$
something similar will hold (up to a small error when 
you transfer the problem from ${\mathbb Z}_N$ to 
${\mathbb Z}$) so that $S_N(1,2,...,d)$ is 
approximately the same as $S_N(M,2M, ..., dM)$.

\item If it's {\it not} possible to bound $m$ in terms of $d$ as asked by the previous problems, give a good construction
of a set of places $x_1, ..., x_d$ such that the number
of corners of $S_N(x_1,...,x_d)$ and $T_N(x_1,...,x_d)$
are rather large.

\item Once the polytope associated with some sequence 
$x_1,...,x_d$ is pinned down we can ask about the 
{\it distribution} of the number of sets $A$ such that
$N^{-1}(1_A*1_A(x_1), ..., 1_A*1_A(x_d)) = 
(t_1,...,t_d)$ for some targets $t_1,...,t_d$.  Perhaps 
it is roughly some kind of multi-dimensional normal
distribution, not just near the mean value but even 
near the boundary of the region.

\end{enumerate}

\section{Proof}

\subsection{Initial setup, defining the weighted graph $G_A$}

We begin by letting $G = (V,E)$ denote the directed de Bruijn graph (with loops) on $2^d$ vertices $V$ labeled by all the binary strings on $d$ symbols, where there is an edge from $v \to w$ if one can produce
the string $w$ by removing the left-most symbol from $v$ and then 
concatenating an additional symbol to right of the string.  
For example,
there is a connection from $10110 \to 01101$, since upon removing
that $1$ on the left of $10110$ we get the string $0110$; and then
adding a $1$ to the right, we get $01101$.

Note that the vertices labeled $000...0$ and $111...1$ both have loops, and are in fact the only vertices with loops.

As is well-known, every vertex $v$ except for $000...0$ and $111...1$
have exactly two edges that point from $v$ to some other vertex; and 
then there are two edges pointing {\it into} $v$ from some other 
vertex.  The vertices $000...0$ and $111...1$ have only one 
edge pointing out to some other vertex and one edge pointing in.

The graph $G$ will be useful for understanding
$S_N(1,2,...,d)$.  To understand $T_N(1,2,...,d)$ we will
need to define the following related graph:  we let 
$G' = (V',E')$ denote the directed graph with vertex labels given by all ordered pairs of binary strings of length $d$ (or alternatively
binary strings of length $2d$) where a vertex with 
label $(s,t)$ ($s$ and $t$ are binary strings 
of length $d$) has an edge pointing to $(s',t')$ if 
there are edges $s' \to s$ and $t \to t'$ in the de Bruijn
graph $G$.  In other words, 
one can produce $s'$ by appending a $0$ or $1$ 
to the beginning (left end) of $s$ while also deleting the right-most 
character from $s$; and one can produce $t'$ by appending
a $0$ or $1$ to the right end of $t$ while also deleting the
left-most character from $t$.

Now we suppose $A \subseteq [N]$ and $A' \subseteq [-N,N]
\cap {\mathbb Z}$.  Note that these gives rise to a vector $\vec h \in H$ via
\begin{equation}\label{hdef}
\vec h\ =\ N^{-1} (1_A*1_{-A}(1), ..., 1_A*1_{-A}(d)).
\end{equation}
and a vector
\begin{equation}\label{h'def}
\vec h'\ =\ (2N+1)^{-1} (1_{A'}*1_{A'}(1), ..., 1_{A'}*1_{A'}(d)).
\end{equation}
We can represent $A$ as a binary
string $s_A$ of length $N$, where there is a $1$ in the $i$th position if
$i \in A$ and a $0$ in the $i$th position if $i \not \in A$.
And we likewise can represent $A'$ by a binary string $s'_{A'}$
of length $2N+1$, indexed by the integers $i \in [-N,N]$,
where there is a $1$ in position $i$ if the element 
$i\in A'$ and there is a $0$ in position $i$ otherwise.

We now map $A$ and $A'$ to weighted directed graphs 
$G_A$ and $G'_{A'}$, respectively, as follows.  
The vertices and edges of $G_A$ are the same as for $G$; and the vertices and edges of $G'_{A'}$ are the same as for $G'$.  It remains to 
say what the edge weights are:  we will first produce from $A$ a walk 
$v_1,v_2,...,v_{N-d+1}$.  We begin by letting $v_1$ be the vertex
whose label is given by the first $d$ symbols in the string
$s_A$.  $v_2$ is then the vertex whose label is given by the 
symbols in positions $2$ through $d+1$ in $s_A$; 
and so on, where $v_{N-d+1}$ is the vertex corresponding to the 
symbols in the $N-d+1$ through $N$ position.

We similarly produce from $A'$ a walk 
$v'_1, ..., v'_{N-d+1}$ through $G'$.  Recall that the 
vertices of $G'$ have labels of the form $(s,t)$, where $s$
and $t$ are binary strings of length $d$.  In our case we
will let $v'_1$ be the vertex where the corresponding string
$s$ is the symbols of $s'_{A'}$ in the positions 
$i=-d+1,...,-2,-1,0$ and the string $t$ is the symbols of 
$s'_{A'}$ in the positions $i=1,2,...,d$.  Then, we let 
$v'_2$ correspond to having $s$ be the symbols in 
positions $i=-d, ..., -1$ in $s'_{A'}$, and $t$ be the 
symbols in positions $i=2,...,d+1$.  We continue in this 
vein, producing $v'_1, v'_2, ..., v'_{N-d+1}$.  
The string $s$ corresponding to $v'_{N-d+1}$ will be
in positions $i = -N+1, ..., -N+d$ and the $t$ will 
correspond to positions $i=N-d+1, ..., N$.  It would seem
that we are missing the contribution of $i=-N$, however
the value of $1_{A'}*1_{A'}(1), ..., 1_{A'}*1_{A'}(d)$
is not affected at all by adding or removing $-N$ from 
the set $A'$, because in order to add to $1,2,...,d$
there $A'$ would also have to contain elements in the
interval $[N+1,...,N+d]$, which is outside its range.

Note that for $i=1,2,...,N-d$ there is an edge $v_i \to v_{i+1}$ in $G$ and likewise an edge
from $v'_i \to v'_{i+1}$ in $G'$.  However, the walks do not necessarily return at the end to $v_1$ or $v'_1$; that is, it isn't necessarily the case that $v_{N-d+1} = v_1$ or
$v'_{N-d+1} = v'_1$.  Nonetheless, by extending the walks to at most $d$ additional vertices $v_{N-d+2}, 
v_{N-d+3}, ...$ and $v'_{N-d+2}, v'_{N-d+3}, ...$
we can produce walk $v_1, v_2, ...,v_M$ and 
$v'_1, v'_2, ..., v'_{M'}$, 
$N-d+1 \leq M, M' \leq N+1$ where $v_M = v_1$ and $v'_{M'} = v'_1$.  

We now count up the number of times during the walk 
$v_1,...,v_M$ that we cross
any particular edge $e \in E$, say this count is $w_e$.
Likewise, for any edge $e' \in E'$ we let 
$w'_{e'}$ denote the number of times during the walk
$v'_1, ..., v'_{M'}$ we cross the edge $e'$.  Then we 
simply let $w_e$ be the weight for the edge $e$ and 
let $w'_{e'}$ be the weight for the edge $e'$. 

Now, for each vertex $v \in V$ the sum of the weights of the edges leading into $v$ from another vertex equals the sum of the weights of edges exiting $v$ to another vertex (note that this excludes the contribution of loops).  This is an immediate 
consequence of the fact that the walk is a closed loop ($v_M = v_1$),
which guarantees that we can pair up each time we enter a vertex $v$ along an edge with a time when we leave it, including for start vertex $v_1$.

Likewise, the sum of the weights of the edges entering any $v' \in V'$ equals the sum of the weights of the edges leaving that vertex.

\subsection{Construction of weighted graphs $G_0, G_1, ...$} \label{weightedgraphs}

We claim that we can decompose $G_A$ into a sum of cycles in the
following sense:  let $C = \{c_1,c_2,...,c_m\}$ denote the set of cycles (including the two loops) in the de Bruijn
graph $G$.  For each $c_i \in C$, we write or identify 
$c_i$ with a formal sum of its
directed edges $c_i = e_{i,1} + e_{i,2} + \cdots + e_{i,\ell_i}$,
where $\ell_i$ is the length of $c_i$ and where $e_{i,j}$ denotes
the $j$th edge of $c_i$.  Likewise, through
an abuse of notation we can identify $G_A$ with a weighted 
formal sum of its edges 
$$
G_A\ =\ \sum_{e \in E} w_e e.
$$
We claim that we can write this as
\begin{equation}\label{claim}
G_A\ =\ \sum_{i=1}^m n_ic_i\ =\ 
\sum_{i=1}^m n_i (e_{i,1} + \cdots + e_{i,\ell_i})
\ =\ \sum_{e \in E}
e \sum_{1\leq i \leq m \atop e\ {\rm an\ edge\ of\ c_i}} n_i,
\end{equation}
where the integers $n_i \geq 0$ satisfy
\begin{equation}\label{Msum}
\sum_{i=1}^m n_i \ell_i\ =\ M-1,
\end{equation}
which is the number of edges in the walk $v_1,...,v_M$.

Likewise, we claim that we can decompose $G'_{A'}$ into
a similar formal sum:
\begin{equation}\label{claim2}
G'_{A'}\ =\ \sum_{i=1}^{m'} n'_i c'_i,
\end{equation}
where $C' = \{c'_1, ..., c'_{m'}\}$ is the set of cycles
in the graph $G'$, where we can write the cycle
$c'_i = e'_{i,1} + \cdots + e'_{i,\ell'_i}$, a formal
sum of its edges.  
Letting $\ell'_i$ denote the length
of the cycle $c'_i$, we also have
\begin{equation}\label{M'sum}
\sum_{i=1}^{m'} n'_i \ell'_i\ =\ M'-1.
\end{equation}

These results (really just a lemma) on decomposing $G_A$ and $G'_{A'}$ into
cycles can be related to the Birkoff-von Neumann theorem \cite{birkhoff},
which says the following.

\begin{theorem}[Birkhoff-von Neumann]  Suppose that 
$Q$ is an $n \times n$ doubly-stochastic matrix, which means that the entries are all in $[0,1]$ and the sum across every row and down any column is always $1$.
Then, $Q$ is a convex combination of $n \times n$
permutation matrices.  That is, there exist permutation
matrices $P_1, ..., P_k$ so that 
$$
Q\ =\ \lambda_1 P_1 + \cdots + \lambda_k P_k,
$$
where all $\lambda_i \geq 0$ and $\lambda_1 + \cdots + 
\lambda_k = 1$.
\end{theorem}

The connection between this theorem and the decompositions
(\ref{claim}) and (\ref{claim2}) is as follows.  First,
we suppose the vertices of $G_A$ are ordered so that
we can talk about the ``$i$th vertex" of the graph.  We
then form the weighted adjacency matrix for $G_A$, where
the $i,j$ entry is the weight of the edge connecting the
$i$th vertex of $G_A$ to the $j$th vertex of $G_A$.  
If $ij$ is not an edge in the graph, the $i,j$ entry
of the matrix is $0$.  Now, the fact that $G_A$ has
the property that the sum of the weights along
edges entering any vertex $v \in V$ is the same as 
the sum of the weights of edges exiting $v$ implies
that the sum of the entries in the $i$th row of the 
adjacency matrix equals the sum of the entries in the
$i$th column of the matrix.  This is not a doubly
stochastic matrix, nor can we easily transform it into
one (simply rescaling rows and then rescaling columns does not work as one might think).  If it {\it were} possible to renormalize in this way to apply
the theorem, writing this matrix as the linear combination of permutation matrices would be the kind of conclusion
we are after.  Each such permutation matrix would
correspond to a disjoint union of cycles (no vertex
appearing in more than one cycle) in the graph.

Instead of applying this theorem, and especially because
the renormalization idea doesn't work in general, we
will just directly prove what we need (the proof is very
simple).  To show that it is possible to write $G_A$ and $G'_{A'}$ this way, we construct sequences 
$$
G_0 := G_A,\ G_1,\ G_2,\ ...\ {\rm and\ } 
G'_0 := G'_{A'},\ G'_1,\ G'_2,\ ...
$$
of weighted graphs by successively removing cycles where all
the edges in the cycle are assigned weight $1$.  In
other words, $G_{i+1} = G_i - c$ for some cycle $c$
of $G_i$, and $G'_{i+1} = G'_i - c'$ for some cycle
$c'$ of $G'_i$.  The fact that such cycles $c$ and $c'$
even exist is due to the fact that at each step starting
at $G_0$ and $G'_0$ we have that the sum of the weights
of edges leading into each vertex of those graphs equals
the sum of the weights leaving.  So at each step $G_i$ 
is either a collection of isolated points (no edges) or else contains a cycle, which we could then subtract away
to get the next graph in the sequence; and the same
for $G'_i$.  Eventually, though, we end up with a 
graph of isolated points, no
edges.  Now for each cycle $c_i$ in $G$ if we let 
$n_i$ denote the number of times we subtracted the cycle
$c_i$ until we reach some $G_i$ to have no edges, then 
we clearly get (\ref{claim}).

Likewise the same argument gives us that (\ref{claim2})
holds for the graph $G'_{A'}$.

\subsection{Building a set from a sum of weighted cycles}
\label{weighted_section}

We will now see how to associate to $\sum_i n_i c_i$ a
special walk that may be different from $v_1,...,v_M$
we had before, and we will do the analogous thing
for $\sum_i n'_i c'_i$ and $v'_1, ..., v'_{M'}$:  
we will think of 
the term $n_1 c_1 = n_1 (e_{1,1} + \cdots + 
e_{1,\ell_1})$ as corresponding to a walk where we start at any vertex 
of $c_1$ and then traverse through the cycle $n_1$ times in total. 
Thus, so far we have a walk of length $n_1 \ell_1$.  Next, we enlarge
the walk by at most $d$ additional steps until we reach any vertex
of $c_2$.  The initial walk on at most $n_1 \ell_1 + d$ edges leading up to a vertex of $c_2$ we will call $W_1$.  Next, we traverse the cycle $c_2$ a
total of $n_2$ times, and then walk to at most $d$ more vertices to 
reach a vertex of $c_3$.  This second walk of length at most 
$n_2 \ell_2 + d$ we will call $W_2$.  At this point we have a walk 
$W_1, W_2$ of total length at most
$n_1 \ell_1 + n_2 \ell_2 + 2d$.  We continue this process, creating walks $W_3, W_4, ..., W_{m-1}$.  When we get to $W_m$ we do not need to
add $d$ additional edges to the end, so it has length $n_m \ell_m$. 
When the process finishes we get a walk $W$ produced
by connecting the smaller walks $W_1, W_2, ..., W_m$ 
end-to-end.  This walk $W$ will have total length at most
\begin{equation}\label{total_walk}
n_1 \ell_1 + \cdots + n_m \ell_m + (m-1) d.
\end{equation}
And the length of $W$ is {\it at least} 
$n_1\ell_1 + \cdots +
n_m \ell_m$, which is the total sum of all 
the weights of $G_A$, which has size $M-1$, where 
$N-d+1 \leq M \leq N+1$.  Note that this implies the 
upper bound in (\ref{total_walk}) satisfies
$$
n_1 \ell_1 + \cdots + n_m \ell_m + (m-1) d\ =\ M-1 + (m-1)d.
$$

We get the same story for $v'_1,...,v'_{M'}$:  we can
create a sequence of walks $W'_1, W'_2, ..., W'_{m'}$,
and then link them together to get one giant walk $W'$ 
of total length at most
\begin{equation} \label{total_walk'}
n'_1 \ell'_1 + \cdots + n'_m \ell'_m + (m-1) d.
\end{equation}
\bigskip

Associated to the walk $W$, we get a corresponding binary string of length at most $M-1 + m d$, which then corresponds to a set $B \subseteq \{1,2,3,...,
M-1+md\}$.  Now, if we let $E$ denote the edges of 
the de Bruijn graph, and let $w(e)$ denote the weight of the
edge in $W$ and $w_0(e)$ denote the weight of the edge
in $G_A$, then we will have for every $e \in E$,
$$
w_0(e)\ \leq\ w(e).
$$
Likewise, associated to $W'$ we get a string 
$B' \subseteq [-(M-1) - md, (M-1) + md]$, 
where the association works in the same way as
when we related a subset $A' \subseteq [-N+d-1, N-d+1]$
to the walk $v'_1,...,v'_{N-d+1}$.  

Let $w'(e')$ denote the
weight of $e'$ when traversing $W'$, and let 
$w'_0(e)$ denote the weight associated to $G'(A')$.  
We will then also have for all such $e'$,
$$
w'_0(e')\ \leq\ w'(e').
$$

Now, each time we move to 
a new vertex in the walk $W$, we either
get a contribution of $0$ or $1$ to $1_B*1_{-B}(j)$; and then adding
up all the contributions of all the edges, as well as the contribution to the convolution coming from the initial choice of vertex, gives us the value of that convolution.  
We get the analogous thing for $1_A*1_{-A}(j)$ traversing the walk used to build $G_A$.

More precisely, we get a contribution of $1$ to $1_B*1_{-B}(j)$ when we walk from a vertex
$v$ to a vertex $v'$ in $W$ if and only if both the following hold:  
(1) that edge $e = vv'$ in 
the de Bruijn graph corresponds to adding $1$ to the right-hand end 
of a string; and, (2) the label (binary string of length $d$) corresponding to $v$ has a $1$ in the $d-j+1$ position.  
And the analogous thing holds for $1_A*1_{-A}(j)$ and the walk used to produce $G_A$.

But this doesn't account for all the possible contributions
to the convolutions $1_B*1_{-B}$ and $1_A*1_{-A}$.  
The remaining possible contribution comes from 
the label of the initial vertex in the walk.  It is equal to the number of pairs of $1$'s that are $j$ apart in the label of the inital vertex.  Thus, this ``remaining possible contribution" would have size at most
$O(d)$.

In the end, $1_A*1_{-A}(j)$ and $1_B*1_{-B}(j)$ are
completely determined by the choice of starting vertex
in the corresponding walks, as well as how often we
visit various edges in the de Bruijn graph.  

We have therefore that for all $j=1,...,d$,
\begin{eqnarray} \label{AB}
1_A*1_{-A}(j)\ &\leq&\ 1_B*1_{-B}(j) + O(d)\nonumber \\
&\leq&\ 1_A*1_{-A}(j) + \sum_{e \in E} (w(e) - w_0(e)) + O(d) \nonumber \\
&\leq&\ 1_A*1_{-A}(j) + (m-1)d + O(d).
\end{eqnarray}

The first inequality follows from the fact that $1_A*1_{-A}(j)$ and $1_B*1_{-B}(j)$ are completely determined by the starting vertex $v_1$ in the walks associated to $A$ and $B$, as well as the number of times each of the edges are visited, where each edge is visited at least as many times in the walk $W$ as in the walk $v_1,...,v_{N-d+1}$.  And the second inequality is due to the fact that for each edge $e$, the difference in $w(e) - w_0(e) \geq 0$ is 
an upper bound for the additional contribution of the edge $e$ to $1_B*1_{-B}(j)$ versus $1_A*1_{-A}(j)$.  

We get a similar story when considering $1_{B'}*1_{B'}(j)$:
to find this convolution we merely need to add up 
the contributions from the starting vertex and then each each edge we traverse in the double de Bruijn graph in completing the walk $W'$.  The contribution of an
edge $(s,t) \to (s',t')$ to the convolution 
will be either $0, 1, $ or $2$.  
Note that the transition $s \to s'$
corresponds to adding or not (depending on whether the
corresponding edge in the de Bruijn graph has label $0$ or $1$) some number 
$1 \leq k \leq M-1 + md$, and the transition $t \to t'$
corresponds to adding or not adding $-k+1$.  

We would get
a contribution of $2$ to $1_{B'}*1_{B'}(j)$ for the
edge $(s,t) \to (s',t')$ if we include $k$ in $B'$ and 
$-k+j$ was already in $B'$, and whether $-k+j$ was included can be read off from the vertex $t$, because its label keeps a record of the previous (in the range
$-k+1$ to $-k+d+1$) several negative elements added to $B'$.  

And traversing $(s,t) \to (s',t')$ gives a contribution of $0$ to $1_{B'}*1_{B'}(j)$ if $k$ is not added to $B'$.  The contribution of $1$ is a special case and corresponds to sums $2b = b+b = j$ for $b \in B'$, so happens only for $1 \leq b \leq j/2$; and so it contributes at most $O(d)$ to the overall value of $1_{B'}*1_{B'}(j)$.  

In the end we get for $j=1,...,d$,
\begin{eqnarray*}
1_{A'}*1_{A'}(j)\ &\leq&\ 1_{B'}*1_{B'}(j) + O(d)\\
&\leq&\ 1_{A'}*1_{A'}(j) + 2 \sum_{e \in E} (w'(e) - w'_0(e)) + O(d) \\
&\leq&\ 1_{A'}*1_{A'}(j) + 2(m-1)d + O(d).
\end{eqnarray*}

\subsection{Relating $B$ to a polytope, proof of
(\ref{AC})} \label{proofac}

Like how we related the walk $W_1,...,W_m$ to the set $B$, we
can individually relate $W_i$ to a set 
$B_i \subseteq \{1,2,..., n_i \ell_i + d\}$, so that for $j=1,...,d$,
\begin{equation}\label{BBj}
1_B*1_{-B}(j)\ =\ \sum_{i=1}^m 1_{B_i}*1_{-B_i}(j) + O(dm).
\end{equation}
Since $W_i$ is a walk around a cycle $c_i$ again and again,
except at the end (where we add vertices to the walk
to hop to a vertex of $c_{i+1}$)
we will have that for $1 \leq b \leq (n_i-1) \ell_i$
that $b \in B_i$ if and only if $b+\ell_i \in B_i$; that is, $B_i$ has a kind of periodicity property.  

We now let, for $i=1,...,m$, $C_i \subseteq {\mathbb Z}_{\ell_i}$ be $B_i \cap [1, \ell_i]$, 
interpreted as a subset mod $\ell_i$.  

From the periodicity property of the set $B_i$ we have
for $j=1,2,...,d$ that
$$
1_{B_i}*1_{-B_i}(j)\ =\ n_i 1_{C_i}*1_{-C_i}(j) + O(d).
$$
Let
\begin{equation}\label{yiell}
\vec y_i\ =\ \ell_i^{-1} (1_{C_i}*1_{-C_i}(1), ..., 
1_{C_i}*1_{-C_i}(d))\ \in\ [0,1]^d.
\end{equation}
For $j=1,...,d$ let
$$
\lambda_j\ :=\ {n_j \ell_j \over M-1}.
$$
Then from (\ref{Msum}) we have 
$\lambda_1 + \cdots + \lambda_m = 1$, and we have 
from (\ref{BBj}) and (\ref{AB}) that for 
any $k=1,2,...,d,$ the $k$th coordinate of 
$\sum_{i=1}^m \lambda_i \vec{y_i}$ equals
$$
1_A*1_{-A}(k)\ =\ 1_B*1_{-B}(k) + O(d)
\ =\ \sum_{i=1}^m n_i 1_{C_i}*1_{-C_i}(k)\ +\ O(dm).
$$
Since $M = N + O(d)$ from (\ref{hdef}) we then deduce that
$$
\vec h\ =\ \sum_{i=1}^m \lambda_i \vec y_i\ +\ O(dm/N).
$$

Thus, we see that the first part of 
(\ref{AC}) holds, and note that the
vectors $\vec y_i$ have rational coordinates as required
by the theorem.

As for the analogous result for $W'_i$, we associate
to it a {\it pair} of sets 
$$
B'_i\ \subseteq\ \{1,2,..., n'_i\ell'_i + d\},
\ B''_i\ \subseteq\ \{-n'_i\ell'_i-d+1, ..., 0\}
$$
so that 
$$
1_{B'}*1_{B'}(j)\ =\ 2\sum_{i=1}^{m'} 1_{B'_i}*1_{B''_i}(j)
+ O(m'd).
$$
We could combine these two sets $B'_i$ and $B''_i$ together into a single set $D_i := B'_i \cup B''_i$ and then write
$$
1_{B'}*1_{B'}(j)\ =\ \sum_{i=1}^{m'} 1_{D_i}*1_{D_i}(j) + O(m'd),
$$
however in order to relate these convolutions $1_{D_i}*1_{D_i}$ to some convolution in a finite
group ${\mathbb Z}_{\ell'_i}$ we need to {\it not} try to
combine $B'_i$ and $B''_i$ into a single set $D_i$.  

Now, as with the set $B_i$, the sets $B'_i$
and $B''_i$ have an approximate periodicity property.
Specifically, for every $1 \leq b \leq (n'_i-1)\ell'_i$
we have that $b \in B'_i$ if and only if 
$b+\ell'_i \in B'_i$; and for 
$-(n'_i-1)\ell' \leq b \leq 0$ we have $b \in B''_i$
if and only if $b-\ell'_i \in B''_i$.

So, as with how we created the sets $C_i$, we let
$$
C'_i\ :=\ B'_i\ \cap\ \{1,2,...,\ell'_i\},\ 
C''_i\ :=\ B''_i\ \cap\ \{-\ell'_i+1, -\ell'_i+2, ..., 0\},
$$
where we are to think of both of these sets as subsets of
${\mathbb Z}_{\ell'_i}$ (instead of just ${\mathbb Z}$).  
Then we observe that
$$
1_{B'_i}*1_{B'_i}(j)\ =\ 2 n'_i 1_{C'_i}*1_{C''_i}(j)
+ O(d).
$$
Let
$$
\vec z_i\ :=\ (\ell'_i)^{-1} (1_{C'_i}*1_{C''_i}(1), ..., 
1_{C'_i}*1_{C''_i}(d))\ \in [0,1]^d.
$$
For $k=1,...,d$ we note that
$$
1_{A'}*1_{A'}(k)\ =\ 1_{B'}*1_{B'}(k) + O(dm') \ =\ 2 \sum_{i=1}^{m'} n'_i 
1_{C'_i}*1_{C''_i}(k)\ +\ O(dm').
$$
So, since $M' = N + O(d)$, from (\ref{h'def})
we deduce that if we let for $i=1,...,m'$,
$$
\lambda'_i\ =\ {n'_i \ell'_i \over M'-1},
$$
then $\lambda'_1 + \cdots + \lambda'_{m'} = 1$ 
and
$$
\vec h'\ =\ \sum_{i=1}^{m'} \lambda'_i \vec z_i
\ +\ O(dm'/N),
$$
which establishes the second part of (\ref{AC}).

\subsubsection{The ``conversely" part of (\ref{AC})}

Now suppose $\vec h \in H$ and $\vec h' \in H'$.  
Thus, there exist $\lambda_1,...,\lambda_m \geq 0$
and $\lambda'_1, ..., \lambda'_{m'} \geq 0$ such that
$$
\lambda_1 + \cdots + \lambda_m\ =\ 1\ =\ 
\lambda'_1 + \cdots + \lambda'_{m'}
$$
and such that
$$
\lambda_1 \vec y_1 \cdots + \lambda_m \vec y_m\ =\ \vec h,
\ {\rm and\ } \lambda'_1 \vec z_1 + \cdots + \lambda'_{m'}
\vec z_{m'}\ =\ \vec h'.
$$

Next, for each $i=1,...,m$ and $j=1,...,m'$ let 
$$
n_i\ :=\ [N \ell_i^{-1} \lambda_i],\ {\rm and\ } 
n'_j\ :=\ [N {\ell'}_j^{-1} \lambda'_j].
$$
And then we consider the weighted de Bruijn graph $G_w$
and the weighted double de Bruijn graph $G'_{w'}$
which we define through linear combinations of cycles
like in section \ref{weightedgraphs} as follows:
$$
G_w\ :=\ \sum_{i=1}^m n_i c_i,\ G'_{w'}\ :=\ \sum_{i=1}^{m'} n'_i c'_i,
$$
where $c_1,...,c_m$ are the cycles of the de Bruijn graph
$G$ and $c'_1, ..., c'_{m'}$ are the cycles of the
double de Bruijn graph $G'$.

Next, we basically repeat the construction of the sequence of walks $W_1,W_2, ..., W_m$ associated to $G_w$ from
the section \ref{weighted_section}, and
then to a subset $B \subseteq \{1,..., M-1 + (m-1)d\}$,
where 
$$
M-1\ =\ \sum_{i=1}^m n_i \ell_i\ =\ N \sum_{i=1}^m 
\lambda_i\ +\ E\ =\ N + E,
$$
where
$$
|E|\ \leq\ \sum_{i=1}^m \ell_i\ \leq\ 2^d m.
$$
Thus, by trimming at most $O(2^d m)$ elements from $B$,
we can ensure that $B \subseteq [1,N]$, and the 
convolutions $1_B*1_{-B}(j)$ will only change by at most
$O(2^d m)$.

Then, for $k=1,2,...,d$, as in section \ref{proofac},
$$
1_B*1_{-B}(k)\ =\ \sum_{i=1}^m n_i 1_{C_i}*1_{-C_i}(k)
\ +\ O(dm).
$$
By (\ref{yiell}), within an error of $O(dm)$ this 
is the $k$th coordinate of $\sum_{i=1}^m \ell_i n_i \vec y_i$, which is, within an error 
$O(m(d + 2^d))$, the $k$th coordinate of
$N \sum_{i=1}^m \lambda_i \vec y_i$.  From this the
first part of (\ref{AC}) follows.

Analogously, to prove the second part of (\ref{AC}) 
we pass from $G'_{w'}$ to walks $W'_1,...,W'_{m'}$
as discussed in section \ref{weighted_section};
and then from these walks we pass to a set
$B' \subseteq [-N, N]$.  And then all the steps above
that worked for $B$ will also work for $B'$, except
that in place of $1_{C_i}*1_{-C_i}(k)$ we have
$2 \cdot 1_{C'_i}*1_{C''_i}(k)$.  In the end, though,
the second part of (\ref{AC}) will follow.

\subsection{Upper bounds on $m$ and $m'$}

One way to bound $m$ and $m'$ would be to bound the
number of cycles in a certain de Bruijn graph and a related graph for $m'$.  However, this will give bounds that are much too large.  

An alternative approach would be to attempt to find some
minimal decomposition of $G_A$ (and $G'_{A'}$) as 
an positive integer linear combination of cycles,
where the number of cycles is minimal.  Such decompositions
might involve significantly fewer cycles than exist
in the de Bruijn graph.  Indeed, it is known
\cite{kulkarni, marcus} that an
$n \times n$ doubly-stochastic matrix can be written
as a sum of {\it at most} $n^2 - 2n+2$ permutation matrices.  However, this bound would apply only for a 
single matrix, not the set of {\it all} doubly-stochastic
matrices at the same time using the same set of $n^2 - 2n + 2$ matrices.

Yet another alternative, which is the one we will actually use, relies on the fact that we don't really need to do something like (related to) bound the number of vertices in the polytope of doubly stochastic matrices.  All we care about is convolutions, and so we can simply use the fact that the $\vec y_i$, 
$i=1,...,m$ and 
$\vec z_j$, $j=1,...,m'$, have rational coordinates
with denominators of size $\ell_i$ and $\ell'_j$,
respectively, where the numerators are integers in 
$\{0,1,...,\ell_i-1\}$ and $\{0,1,...,\ell'_i-1\}$, respectively.  So, an upper bound for the number of 
vectors $\vec x_1, ..., \vec x_m$ is 
$$
\sum_{\ell_i} \ell_i^d,
$$
where the sum is over all the $\ell_i$ that are possible
cycle lengths in a de Bruijn graph.  

As is well known, de Bruijn graphs contain Hamilton cycles, so we do not get any better bound than $\ell_i \leq 2^d$
on the possibilities for the length $\ell_i$.  Thus, we
get the upper bound
$$
m\ \leq\ \sum_{j=1}^{2^d} j^d\ \leq\ 2^{d(d+1)}.
$$
We can actually improve this by a factor of $d$ or so, 
but there is no reason to bother since the bound is probably
nowhere near the true upper bound.

Using an analogous argument and the fact that 
the double de Bruijn graph has $4^d$ vertices 
we get that
$$
m'\ \leq\ \sum_{j=1}^{4^d} j^d\ \leq\ 4^{d(d+1)}.
$$


\begin{thebibliography}{999}

\bibitem{birkhoff} G. Birkhoff, Tres observaciones sobre el algebra lineal, Univ. Nac. Tucumán, Rev. Ser. A, no. 5 (1946), 147–151. 

\bibitem{bourgain} J. Bourgain, {\it On arithmetic progressions in sums of sets of integers}, A tribute to Paul Erd\H os, 105–109 (CUP, 1990)

\bibitem{croot1} E. Croot, I. Ruzsa, and T. Schoen, 
{Arithmetic progressions in sparse sumsets}, Combinatorial
Number Theory, 157-164 (de Gruyter, Berlin, 2007).

\bibitem{croot2} E. Croot and O. Sisask, 
{\it A probabilistic technique for finding almost-periods of convolutions}, Geom. Funct. Anal. {\bf 20} (2010), no. 6, 1367-1396.

\bibitem{green} B. Green, {\it Arithmetic progressions in sumsets}, Geom. Funct. Anal. 12 (2002), no. 3, 584–597.

\bibitem{kulkarni} J. Kulkarni, E. Lee, and M. Singh, Minimum Birkhoff-von Neumann Decomposition, In: Eisenbrand, F., Koenemann, J. (eds) Integer Programming and Combinatorial Optimization. IPCO 2017. Lecture Notes in Computer Science, vol 10328.

\bibitem{marcus} M. Marcus, R. Ree, {\it Diagonals of doubly stochastic matrices}, Q. J. Math. {\bf 10} (1959), 296–302.


\end{thebibliography}
\end{document}